\documentclass[11pt]{amsart}
\usepackage[english]{babel}
\usepackage[utf8x]{inputenc}
\usepackage{amsmath,amsfonts,stackrel,amsthm,enumerate}
\usepackage{graphicx,setspace,mathrsfs}
\usepackage[margin=1in]{geometry}
\usepackage{xcolor}
\usepackage{amssymb,lineno}
\newtheorem{definition}{Definition}
\newtheorem{theorem}{\bf Theorem}[section]
\newtheorem{remark}{\bf Remark}[section]
\newtheorem{proposition}{Proposition}[section]

\newtheorem{corollary}{Corollary}[section]

\newcommand{\process}{{\normalfont\textsc{GSFPP-VO}}}
\newcommand{\processf
}{{Generalized Space-Fractional Poisson Process via Variable-Order Stable Subordinator}}
\newcommand{\gsfppvo}{\{\xi^{\alpha(t)}(t)\}_{t\geq0}}
\usepackage[colorinlistoftodos]{todonotes}
\usepackage{enumitem}
\usepackage{enumerate}
\title{\processf
}
\author{Reetendra Singh}
\address{Decision Science Area, Indian Institute of Management Visakhapatnam}
\email{reetendra.singh21-02@iimv.ac.in}
\author{Aditya Maheshwari}
\address{Operations Management \& Quantitative Techniques Area, Indian Institute of Management Indore}
\email{adityam@iimidr.ac.in}

\begin{document}

\begin{abstract}
  This paper introduces a variable-order stable subordinator (VOSS)  $S^{\alpha(t)}(t)$ with index $\alpha(t)\in(0,1)$, where $\alpha(t)$ is a right-continuous piecewise constant function. We drive the \processf ~(\process) defined by $\{N(S^{\alpha(t)}(t))\}_{t \geq 0}$, obtained by time-changing a homogeneous Poisson process $\{N(t,\lambda)\}_{t\geq 0}$ with rate parameter $\lambda>0$ by an independent VOSS. Explicit expressions for the Laplace transform, probability generating function, probability mass function, and moment generating function of the \process\ are derived, and these quantities are shown to satisfy  partial differential equations. Finally, we establish the associated generalized distributions, analyze the hitting-time properties, and characterize the L\'evy measures of the \process.
\end{abstract}
	\maketitle
 \section{Introduction}
 In recent years, the fractional Poisson process (FPP) has been considered a significant stochastic model for count data with notable mathematical properties. These models have applications across various disciplines, including epidemiology, industry, biology, queuing theory, traffic flow, reliability, finance, economics, statistical physics, and business \cite{haight1967handbook,byrne1969properties,golding2006physical,stanislavsky2008two,bening2012generalized,wang2022explicit,soni2023bivariate}. Laskin \cite{laskin2009some} introduced the fractional Poisson process (FPP) to establish a novel set of quantum coherent states, along with a fractional extension of Bell polynomials, numbers, and Stirling's numbers of the second kind, which was revolutionary work in the field of FPP. Later on, the space-time inclusion with non-homogeneous FPP has been elaborated by Maheshwari et al. \cite{maheshwari2019non}. Furthermore, fractional Poisson processes have been used to analyse renewal processes with inter-times between events that follow Mittag-Leffler distributions, as demonstrated and involves substituting the time derivative in the equations that determine the state probabilities with the fractional derivative, specifically in the sense of Caputo \cite{mainardi04,beghin2009fractional}.

 Orsingher et al. \cite{Orsingher2012} introduced the space-fractional Poisson process $N^{\alpha}(t)$ and demonstrated that this process can be seen as a homogeneous Poisson process $N(t)$, which is subordinated to a positively skewed stable process $S^{\alpha}(t)$ with a Laplace transform. Later on, Polito et al. \cite{polito2016generalization} presented the generalization of the space fractional Poisson process and its connections to some Levy processes without consideration of time in subordinator. The underpinned literature on the space-fractional Poisson process (SFPP) shows a scope of studies concerning the subordinator as a function of time. To fill this gap and add some contributions to the study of Orsingher et al. \cite{Orsingher2012}, we propose generalizing the space-fractional Poisson process with a time-varying subordinator.
 
In this study, we introduce a generalization of the space-fractional Poisson process via variable stable subordinator (\process). This generalization offers intriguing insights into the characteristics and properties of the process. This proposed process satisfies the condition of independent identical properties with infinite moments. Here, we provide a comprehensive outline of the generalization of the space-fractional Poisson process via a variable stable subordinator. We derive the time-varying space-fractional Poisson process and solve Laplace differential equations for the probability mass function (pmf), probability generating function (pgf), and moment generating function (mgf) for the order of subordinator $\alpha(t)$ with time $t$. Then, finally, we investigate the process's hitting time for \process.
 
The structure of this paper is organized as follows: Section \ref{sec:prelims} of the paper contains preliminary notations and definitions of other processes. In Section \ref{sec:tvsfpp_hittingtime}, we obtained a \processf\ with the distributions and differentiation. Finally, we concluded with the examination of hitting time, and some other properties of the \process.

 \section{Preliminaries}\label{sec:prelims}
 \noindent
 In this section, we introduce some of the notations and results that could be used throughout the paper. Let $\mathbb{N}=\{1,2,\ldots,n\}$ be set of natural numbers.\par
 \noindent
 There are some special functions that will be required in the later part of the study.
 A space-fractional Poisson process by means of the fractional difference operator introduced by Garra et al. \cite{garra2017some}
 \begin{align}\label{operator}
     \Delta^{\alpha}=(1-B)^{\alpha}, ~\alpha \in(0,1].
 \end{align}
 This frequently comes up when investigating long-memory time series.
 The state probabilities $p_k^{\alpha}(t)$ are dependent on all probabilities by the operator (\ref{operator}). 
The classical homogeneous Poisson process is recovered for $\alpha=1$, and the state probabilities $p_k(t)$ rely only on $p_{k-1}(t)$. Distribution for the space-fractional Poisson process has been established by 
\begin{align}\label{space_frac}
    p_k^{\alpha}(t)=Pr\{N^{\alpha}(t)=k\}=\frac{(-1)^k}{k!}\sum_{r=0}^{\infty}\frac{(-{\lambda}^{\alpha} r)}{r!} \frac{\Gamma{(\alpha r+1})}{\Gamma{(\alpha r+1-k})},~k\geq 0.
\end{align}

The one parameter Mittag–Leffler function
\begin{align}
    E_{\alpha}(x)=\sum_{r=0}^{\infty} \frac{x^r}{\Gamma(\alpha r+1)}.
\end{align}
Mittag-Leffler function is very important function to find the probability generating function and other distributions.\\
Let $\{S(t)\}_{t\geq 0}$ is a non-decreasing L\'evy process with Laplace transform, which is called L$\acute{e}$vy subordinator \cite{soni2023bivariate}, if
\begin{align*}
    \mathbb{E}(e^{-wS(t)})=e^{-t\phi(w)},~w\geq0,
\end{align*}
where $\phi(w)$ is the Laplace proponent denoted by
\begin{align*}
    \phi(w)=aw+\int_0^{\infty}(1-e^{-ws})\nu(ds),~a\geq0.
\end{align*}
Here, $a$ is the drift coefficient and $\nu$ is a non-negative L$\acute{e}$vy measure on the half line, which satisfies the following inequality
\[
    \int_{0}^{\infty}\min\{s,1\}\nu(ds)<\infty,~ and~ \nu([0,\infty))=\infty.
\]
In a such way that $\{S(t)_{t\geq 0}\}$ almost surely (a.s.) has strictly growing sample paths.

\begin{definition}
	Let $\{N(t,\lambda)\}_{t\geq 0}$ be the homogeneous Poisson process with rate parameter $\lambda>0$ and $\{S^{\alpha}(t)\}_{t\geq 0}$ be an $\alpha-$stable subordinator with index $\alpha\in (0,1)$, independent of $\{N(t,\lambda)\}_{t\geq 0}$. Then, define $ N(S^{\alpha}(t))$ be the space fractional Poisson process with $\alpha$-stable subordinator.
\end{definition}
The section introduce the \processf\ with possible properties.
 \section{Time-varying space fractional Poisson process and Hitting time}\label{sec:tvsfpp_hittingtime}
 \noindent

\begin{definition}[Variable-order stable subordinator]
Let $\alpha: [0,\infty) \to (0,1)$ be a right-continuous, piecewise constant function,  i.e. for any finite time horizon $t > 0$, there exists a partition  $
P = \{0 = t_0 < t_1 < t_2 < \cdots < t_n = t\}$
of the interval $[0,t]$ such that
\begin{equation}\label{order}
\alpha(t) = \alpha_k \quad \text{for all } t \in [t_{k-1}, t_k), \quad k = 1, \dots, n,
\end{equation}
where each $\alpha_k \in (0,1)$ is constant on the $k$-th subinterval.

A stochastic process $\{S^{\alpha(\cdot)}(t)\}_{t \geq 0}$ is called a \emph{variable order stable subordinator (VOSS)
}  if it is constructed as the cumulative sum of independent classical stable subordinators on each subinterval of the (global) partition induced by the jumps of $\alpha(\cdot)$.  In other words, 
for any $t \geq 0$, let $k$ be such that $t_{k-1} \leq t < t_k$ (extending the partition to $[0,\infty)$ as needed), then
\[
S^{\alpha(\cdot)}(t) = \sum_{j=1}^{k-1} S^{\alpha_j}(t_j - t_{j-1}) + S^{\alpha_k}(t_k - t_{k-1}),
\]
where $S^{\alpha_j}(\cdot)$ is an independent classical $\alpha_j$-stable subordinator (i.e., a Lévy process with Laplace transform $\mathbb{E}[e^{-\lambda S^{\alpha_j}(t)}] = e^{-t \lambda^{\alpha_j}}$ for $t \geq 0$, $\lambda \geq 0$). Note that the  time $t$ for both $\alpha$ and $\xi$ are same, which is the variable for stable subordinator.
\end{definition}

It can be noted that the process satisfies:
\begin{itemize}
    \item[(a)] $S^{\alpha(\cdot)}(0) = 0$ almost surely.
    \item[(b)] Independent increments: for any $0 \leq t_0 < t_1 < \cdots < t_m$, the increments $S^{\alpha(\cdot)}(t_1) - S^{\alpha(\cdot)}(t_0), \dots, S^{\alpha(\cdot)}(t_m) - S^{\alpha(\cdot)}(t_{m-1})$ are independent.
    \item[(c)] The distribution of each increment over an interval $[s,r)$ depends only on the piecewise constant values of $\alpha(\cdot)$ restricted to $[s,r)$ and the length of subintervals.
    \item[(d)] For $s < r$ and $\lambda \geq 0$, the Laplace transform is
    \[
    \mathbb{E}\left[\exp\left(-\lambda \left(S^{\alpha(\cdot)}(r) - S^{\alpha(\cdot)}(s)\right)\right)\right] = \exp\left( -\sum_{j: [t_{j-1}, t_j) \cap [s,r) \neq \emptyset} (t_j' - t_{j-1}') \lambda^{\alpha_j} \right),
    \]
    where $[t_{j-1}', t_j') = [t_{j-1}, t_j) \cap [s,r)$ (adjusted for partial overlap at endpoints).
\end{itemize}

\noindent \textbf{Note:} When $\alpha(t) \equiv \alpha$ (constant, corresponding to a trivial partition with $n=1$), the process reduces to the classical $\alpha$-stable subordinator, which is a Lévy process.
 
 \begin{definition} Let ${N(t,\lambda)}_{t\geq 0}$ be Poisson process with rate $\delta$ and $\{S^{\alpha(\cdot)}(t)\}_{t \geq 0}$ be a VOSS with index $\alpha(t)$, which is independent of  ${N(t,\lambda)}_{t\geq 0}$. We define \processf\ (\process) with index $\alpha(t)\in (0,1)$  as time-changed Poisson process subordinated with an independent VOSS, i.e.
 \begin{equation}
      \gsfppvo=N(S^{\alpha(t)}(t)),t\geq 0.
 \end{equation}
\end{definition}

	\begin{proposition}\label{prop_lap}
		Let $\{S^{\alpha(t)}(t)\}_{t\geq 0}$ be the variable-order stable subordinator, where $\alpha(t)$ is function of time as defined in \eqref{order}. Then its Laplace transform is given by
		\begin{equation}
		\mathbb{E}[e^{-\lambda S^{\alpha(t)}(t)}] = e^{-\int_0^{t}\lambda^{\alpha(s)}ds}.
		\end{equation}	
   \end{proposition}
	\begin{proof}
    Using $\alpha$-stable subordinator, which has partition of time interval (see eq \eqref{order})
		\begin{align*}
		\mathbb{E}[e^{-\lambda S^{\alpha(t)}(t)}] =& \mathbb{E}[e^{-\lambda\sum_{k=1}^n[S^{\alpha_k}(t_k)-S^{\alpha_k}(t_{k-1})}]\\ =& \prod_{k=1}^n e^{-\lambda^{\alpha_k}(t_k-t_{k-1})}\\ =& e^{-\sum \lambda^{\alpha_k}(t_k-t_{k-1})}\\ =& e^{-\int_0^t \lambda^{\alpha(s)}ds}.\qedhere
		\end{align*}
	\end{proof}  
	\noindent
	Let $p_k^{\alpha(t)}(t)= \mathbb{P}[\gsfppvo=k]$ be the \textit{pmf} of \process, which satisfies the following differential equation
	\begin{equation}
	\frac{d}{dt}p_k^{\alpha(t)}(t)=  -\lambda^{\alpha(t)}(1-B)^{\alpha(t)}p_k^{\alpha(t)}(t),
	\end{equation} 
	$$ p_k^{\alpha(0)}(0)=\begin{cases}
	1,~~k=0,\\ 0,~~k>0.
	\end{cases}$$

 \begin{equation}
   \begin{cases}
       \frac{d}{dt}p_k^{\alpha(t)}(t)=  -\lambda^{\alpha(t)}\sum_{r=0}^{k}\frac{\Gamma(\alpha(t)+1)}{r!\Gamma(\alpha(t)+1-r)}(-1)^r p_{k-r}^{\alpha(t)}(t),\\
       p_k^{\alpha(0)}(0)=\begin{cases}
	1,~~k=0,\\ 0,~~k>0.
   \end{cases}
   \end{cases}
 \end{equation}

	\begin{theorem}
		Let $\psi_{\alpha(t)}(u,t)$ be the pgf of \process\ and it satisfies the following partial differential equation.
		\begin{equation} \label{pde}
		\frac{\partial}{\partial t}\psi_{\alpha(t)}(u,t) =-\lambda^{\alpha(t)}\psi_{\alpha(t)}(u,t)(1-u)^{\alpha(t)},
		\end{equation}
		$$ \psi_{\alpha(0)}(u,0)=1.$$
	\end{theorem}
	\begin{proof}
		
        We have $\psi_{\alpha(t)}(u,t)= \mathbb{E}[u^{\gsfppvo}]=\sum_{m=0}^{\infty}u^m p_m^{\alpha(t)}(t)$. Then
		\begin{align*}
		\frac{\partial}{\partial t}\psi_{\alpha(t)}(u,t) =& \sum_{m=0}^{\infty}u^m \frac{d}{dt} p_m^{\alpha(t)}(t)\\
		=& -\sum_{m=0}^{\infty}u^m \lambda^{\alpha(t)}(1-B)^{\alpha(t)}p_k^{\alpha(t)}(t)\\
		=& -\sum_{m=0}^{\infty}u^m \lambda^{\alpha(t)}\sum_{r=0}^{m}\frac{(-1)^r}{r!}\frac{\Gamma(\alpha(t)+1)}{\Gamma(\alpha(t)+1-r)}p_{m-r}^{\alpha(t)}(t)\\
		=&-\sum_{r=0}^{\infty}\sum_{m=r}^{\infty}u^m \lambda^{\alpha(t)}\frac{(-1)^r}{r!}\frac{\Gamma(\alpha(t)+1)}{\Gamma(\alpha(t)+1-r)}p_{m-r}^{\alpha(t)}(t)\\
		=&-\sum_{r=0}^{\infty}\sum_{m=0}^{\infty}u^{m+r} \lambda^{\alpha(t)}\frac{(-1)^r}{r!}\frac{\Gamma(\alpha(t)+1)}{\Gamma(\alpha(t)+1-r)}p_{m}^{\alpha(t)}(t)\\
		=&-\lambda^{\alpha(t)}\psi_{\alpha(t)}(u,t)\sum_{r=0}^{\infty}u^r\frac{(-1)^r}{r!}\frac{\Gamma(\alpha(t)+1)}{\Gamma(\alpha(t)+1-r)}\\
		=&-\lambda^{\alpha(t)}\psi_{\alpha(t)}(u,t)(1-u)^{\alpha(t)}.
		\end{align*}
	\end{proof}
	\begin{remark} \label{pgf}
		From equation (\ref{pde}), we observe that
		\begin{align*}
		\psi_{\alpha(t)}(u,t)=& e^{-\int_0^{t}\lambda^{\alpha(s)}(1-u)^{\alpha(s)}ds}\\
		=& e^{-\sum_{k=1}^n\lambda^{\alpha_k}(1-u)^{\alpha_k}(t_k-t_{k-1})}\\
		=& \prod_{k=1}^n e^{-\lambda^{\alpha_k}(1-u)^{\alpha_k}(t_k-t_{k-1})},
		\end{align*}
		which is the \textit{pgf} of $\Gamma(t_1-t_0)+\Gamma(t_2-t_2)+\ldots+\Gamma(t_n-t_{n-1}).$ \\
		Hence $$ \gsfppvo \stackrel{d}{=} \Gamma(t_1-t_0)+\Gamma(t_2-t_2)+\ldots+\Gamma(t_n-t_{n-1}). $$
	\end{remark}
    \noindent
    Next, we discuss the pmf and mgf of the \process.
	\begin{theorem}
		Let $\gsfppvo$ be the \process, then its \textit{pmf} can be expressed as
		\begin{align}
		\mathbb{P}[\gsfppvo=m]=& \frac{(-1)^m}{m!}\sum_{r=0}^{\infty}\frac{(-1)^r}{r!}\sum_{\substack{x_1,x_2,\ldots x_n \geq 0\\ x_1+x_2+\ldots+x_n=r}}\frac{r!}{x_1!x_2!\ldots x_n!}\lambda^{\sum_{i=1}^n\alpha_ix_i}(\alpha_1x_1+\ldots+\alpha_nx_n)_m\nonumber \\&  \prod_{k=1}^n (t_k-t_{k-1})^{x_k}.
		\end{align}
	\end{theorem}
	\begin{proof}
		From the Remark \ref{pgf}, we have that
		\begin{align*}
		\psi_{\alpha(t)}(u,t)=& e^{-\sum_{k=1}^n\lambda^{\alpha_k}(1-u)^{\alpha_k}(t_k-t_{k-1})}\\
		=& \sum_{r=0}^{\infty}\frac{(-1)^r}{r!}\left[\sum_{k=1}^n\lambda^{\alpha_k}(1-u)^{\alpha_k}(t_k-t_{k-1})\right]^r\\
		=&\sum_{r=0}^{\infty}\frac{(-1)^r}{r!}\sum_{\substack{x_1,x_2,\ldots x_n \geq 0\\ x_1+x_2+\ldots+x_n=r}}\frac{r!}{x_1!x_2!\ldots x_n!}[(\lambda(1-u))^{\alpha_1}(t_1-t_0)]^{x_1}\\
		&[(\lambda(1-u))^{\alpha_2}(t_1-t_0)]^{x_2}\ldots[(\lambda(1-u))^{\alpha_n}(t_1-t_0)]^{x_n}\\
		=& \sum_{r=0}^{\infty}\frac{(-1)^r}{r!}\sum_{\substack{x_1,x_2,\ldots x_n \geq 0\\ x_1+x_2+\ldots+x_n=r}}\frac{r!}{x_1!x_2!\ldots x_n!}\lambda^{\sum_{i=1}^n\alpha_ix_i} \\
		&(1-u)^{\sum_{i=1}^n\alpha_ix_i}\prod_{k=1}^n (t_k-t_{k-1})^{x_k} \\
		=&\sum_{r=0}^{\infty}\frac{(-1)^r}{r!}\sum_{\substack{x_1,x_2,\ldots x_n \geq 0\\ x_1+x_2+\ldots+x_n=r}}\frac{r!}{x_1!x_2!\ldots x_n!}\lambda^{\sum_{i=1}^n\alpha_ix_i} \\&
		\sum_{m=0}^{\infty}\frac{(-1)^m}{m!}(\alpha_1x_1+\ldots+\alpha_nx_n)_m u^m \prod_{k=1}^n (t_k-t_{k-1})^{x_k}\\
		=&\sum_{m=0}^{\infty}\frac{(-1)^m}{m!}u^m\sum_{r=0}^{\infty}\frac{(-1)^r}{r!}\sum_{\substack{x_1,x_2,\ldots x_n \geq 0\\ x_1+x_2+\ldots+x_n=r}}\frac{r!}{x_1!x_2!\ldots x_n!}\lambda^{\sum_{i=1}^n\alpha_ix_i}\\
		&(\alpha_1x_1+\ldots+\alpha_nx_n)_m u^m \prod_{k=1}^n (t_k-t_{k-1})^{x_k},
		\end{align*}
        which completes the proof.
	\end{proof}
\begin{theorem}
		Let $\gsfppvo$ be the \process\ then summation of its \textit{pmf} should be equal to one, which can be expressed as
		\begin{equation*}
  \sum_{m=0}^{\infty}p_m^{\alpha(t)}(t)=1.
		\end{equation*}
  \end{theorem}
 
\begin{proof}
By the definition of pmf, we have that
    \begin{align*}
       \sum_{m=0}^{\infty}p_m^{\alpha(t)}(t)=& \sum_{m=0}^{\infty}\mathbb{P}[\gsfppvo=m]\\=&  \sum_{m=0}^{\infty} \frac{(-1)^m}{m!}\sum_{r=0}^{\infty}\frac{(-1)^r}{r!}\sum_{\substack{x_1,x_2,\ldots x_n \geq 0\\ x_1+x_2+\ldots+x_n=r}}\frac{r!}{x_1!x_2!\ldots x_n!}\lambda^{\sum_{i=1}^n\alpha_ix_i}(\alpha_1x_1+\ldots+\alpha_nx_n)_m \\& \prod_{k=1}^n (t_k-t_{k-1})^{x_k}
          \\=& \sum_{m=0}^{\infty} \frac{(-1)^m}{m!}\sum_{r=0}^{\infty}\frac{(-1)^r}{r!}\sum_{\substack{x_1,x_2,\ldots x_n \geq 0\\ x_1+x_2+\ldots+x_n=r}}\frac{r!}{x_1!x_2!\ldots x_n!}(\lambda_{1}^{\alpha_{1}})^{x_1}(\lambda_{2}^{\alpha_{2}})^{x_2}\ldots(\lambda_{n}^{\alpha_{n}})^{x_n}\\&(\alpha_1x_1+\ldots+\alpha_nx_n)_m  (t_1-t_{0})^{x_1}(t_2-t_{1})^{x_2}\ldots (t_n-t_{n-1})^{x_n}\\
          =& \sum_{m=0}^{\infty} \frac{(-1)^m}{m!}\sum_{r=0}^{\infty}\frac{(-1)^r}{r!}\sum_{\substack{x_1,x_2,\ldots x_n \geq 0\\ x_1+x_2+\ldots+x_n=r}}\frac{r!}{x_1!x_2!\ldots x_n!}({\lambda_{1}^{\alpha_{1}}}(t_{1}-t_{0}))^{x_1}({\lambda_{2}^{\alpha_{2}}}(t_{2}-t_{1}))^{x_2}\ldots \\&({\lambda_{n}^{\alpha_{n}}}(t_{n}-t_{n-1}))^{x_n}(\alpha_1x_1+\ldots+\alpha_nx_n)_m \\
           =& \sum_{r=0}^{\infty}\frac{(-1)^r}{r!}\sum_{\substack{x_1,x_2,\ldots x_n \geq 0\\ x_1+x_2+\ldots+x_n=r}}\frac{r!}{x_1!x_2!\ldots x_n!}({\lambda_{1}^{\alpha_{1}}}(t_{1}-t_{0}))^{x_1}({\lambda_{2}^{\alpha_{2}}}(t_{2}-t_{1}))^{x_2}\ldots \\&({\lambda_{n}^{\alpha_{n}}}(t_{n}-t_{n-1}))^{x_n}\sum_{m=0}^{\infty} \frac{(-1)^m}{m!}(\alpha_1x_1+\ldots+\alpha_nx_n)_m \\
           =& \sum_{r=0}^{\infty}\frac{(-1)^r}{r!}\sum_{\substack{x_1,x_2,\ldots x_n \geq 0\\ x_1+x_2+\ldots+x_n=r}}\frac{r!}{x_1!x_2!\ldots x_n!}({\lambda_{1}^{\alpha_{1}}}(t_{1}-t_{0}))^{x_1}({\lambda_{2}^{\alpha_{2}}}(t_{2}-t_{1}))^{x_2}\ldots \\&({\lambda_{n}^{\alpha_{n}}}(t_{n}-t_{n-1}))^{x_n}(1-1)^{(\alpha_1x_1+\ldots+\alpha_nx_n)}.
           \intertext{This implies that $\alpha_1x_1+\ldots+\alpha_nx_n=0\Rightarrow x_1=0=\ldots=x_n$, we have that}
           =& 1.\qedhere
    \end{align*}
\end{proof}

\noindent
\begin{theorem}
 Let $\gsfppvo$ be the \process\ then its \textit{mgf} can be expressed as  
 \begin{align*}
     M_{{\gsfppvo}}(t)= \sum_{r=0}^{\infty} \frac{(-1)^r}{r!}\frac{\Gamma(\sum_{i=1}^{n}x_{i}+1)}{\prod_{i=1}^{n}\Gamma(x_{i}+1)} \prod_{i=1}^{n}(\lambda_{i}^{\alpha_i}(t_i-t_{i-1}))^{x_i} x_{1}^{-\alpha_1}\ldots x_{n}^{-\alpha_n}.
 \end{align*}
\end{theorem}

\begin{proof}
By the definition of mgf, we can write as
    \begin{align*}
     M_{{\gsfppvo}}(t)=&\mathbb{E}[e^{t{{\gsfppvo}}}] \\=& \sum_{m=0}^{\infty} e^{tm} p_m^{\alpha(t)}(t)\\=&
   \sum_{m=0}^{\infty}e^{tm}\frac{(-1)^m}{m!}\sum_{r=0}^{\infty}\frac{(-1)^r}{r!}\sum_{\substack{x_1,x_2,\ldots x_n \geq 0\\ x_1+x_2+\ldots+x_n=r}}\frac{r!}{x_1!x_2!\ldots x_n!}\lambda^{\sum_{i=1}^n\alpha_ix_i}(\alpha_1x_1+\ldots+\alpha_nx_n)_m \\& \prod_{k=1}^n (t_k-t_{k-1})^{x_k}\\=&\sum_{m=0}^{\infty}e^{tm} \frac{(-1)^m}{m!}\sum_{r=0}^{\infty}\frac{(-1)^r}{r!}\sum_{\substack{x_1,x_2,\ldots x_n \geq 0\\ x_1+x_2+\ldots+x_n=r}}\frac{r!}{x_1!x_2!\ldots x_n!}(\lambda_{1}^{\alpha_{1}})^{x_1}(\lambda_{2}^{\alpha_{2}})^{x_2}\ldots(\lambda_{n}^{\alpha_{n}})^{x_n} \\&(\alpha_1x_1+\ldots+\alpha_nx_n)_m (t_1-t_{0})^{x_1}(t_2-t_{1})^{x_2}\ldots (t_n-t_{n-1})^{x_n}\\=& \sum_{m=0}^{\infty}e^{tm} \frac{(-1)^m}{m!}\sum_{r=0}^{\infty}\frac{(-1)^r}{r!}\sum_{\substack{x_1,x_2,\ldots x_n \geq 0\\ x_1+x_2+\ldots+x_n=r}}\frac{r!}{x_1!x_2!\ldots x_n!}({\lambda_{1}^{\alpha_{1}}}(t_{1}-t_{0}))^{x_1}({\lambda_{2}^{\alpha_{2}}}(t_{2}-t_{1}))^{x_2}\ldots \\&({\lambda_{n}^{\alpha_{n}}}(t_{n}-t_{n-1}))^{x_n}(\alpha_1x_1+\ldots+\alpha_nx_n)_m \\=& \sum_{r=0}^{\infty}\frac{(-1)^r}{r!}\sum_{\substack{x_1,x_2,\ldots x_n \geq 0\\ x_1+x_2+\ldots+x_n=r}}\frac{r!}{x_1!x_2!\ldots x_n!}({\lambda_{1}^{\alpha_{1}}}(t_{1}-t_{0}))^{x_1}({\lambda_{2}^{\alpha_{2}}}(t_{2}-t_{1}))^{x_2}\ldots \\&({\lambda_{n}^{\alpha_{n}}}(t_{n}-t_{n-1}))^{x_n}\sum_{m=0}^{\infty}e^{tm} \frac{(-1)^m}{m!}(\alpha_1x_1+\ldots+\alpha_nx_n)_m 
   \\=&
   \sum_{r=0}^{\infty} \frac{(-1)^r}{r!}\frac{\Gamma(\sum_{i=1}^{n}x_{i}+1)}{\prod_{i=1}^{n}\Gamma(x_{i}+1)} \prod_{i=1}^{n}(\lambda_{i}^{\alpha_i}(t_i-t_{i-1}))^{x_i} x_{1}^{-\alpha_1}\ldots x_{n}^{-\alpha_n}.\qedhere
\end{align*}
\end{proof}
\noindent
\begin{corollary}
Let $\gsfppvo$ be the \process\ then its probability characteristic function (\textit{pcf}) is
\begin{align*}
    \phi_{{\gsfppvo}}(t)=&\mathbb{E}[e^{it{{\gsfppvo}}}] \\=& \sum_{m=0}^{\infty} e^{itm} p_m^{\alpha(t)}.
\end{align*}
\end{corollary}
\noindent
Now, we discuss the hitting-time properties of \process.
\begin{definition}
Let $\tau_k$ be a sequence of the first hitting time for the \process\ $\gsfppvo$ and $S^{\alpha(t)}(t)$ is stable subordinator with order of subordinator $\alpha(t)$, which is function of time. The hitting probability $\mathbb{P}[\tau^{\alpha(t)}_k<\infty]$ of the \process\ $\gsfppvo$ is defined as
\begin{align*}
    \mathbb{P}[\tau^{\alpha(t)}_{k}<\infty]:=\mathbb{P}(\inf\{t:\gsfppvo=k\}<\infty).
\end{align*}
\end{definition}
Since, it is given that $\mathbb{P}\{\tau^{\alpha(t)}_{k}(t)<t\}=\mathbb{P}\{\gsfppvo \geq k\}$ and we have the distribution of the \process\ for hitting time 
\begin{align*}
    \mathbb{P}\{\tau^{\alpha(t)}_k<t\}=& \sum_{m=k}^{\infty}\frac{(-1)^m}{m!}\sum_{r=0}^{\infty}\frac{(-1)^r}{r!}\sum_{\substack{x_1,x_2,\ldots x_n \geq 0\\ x_1+x_2+\ldots+x_n=r}}\frac{r!}{x_1!x_2!\ldots x_n!}\lambda^{\sum_{i=1}^n\alpha_ix_i}(\alpha_1x_1+\ldots+\alpha_nx_n)_m\\& \prod_{k=1}^n (t_k-t_{k-1})^{x_k}.
\end{align*}
We know that cumulative distribution function $F_{\tau_{k}^{\alpha(t)}}(t):\mathbb{R}^n\rightarrow[0,1]$ defines as $F_{\tau_{k}^{\alpha(t)}}(t)=\mathbb{P}(\tau_{1}^{\alpha_{t_1}}\leq t_{1},\tau_{2}^{\alpha_{t_2}}\leq t_{2},\ldots,\tau_{k}^{\alpha_{t_k}}\leq t_{k})$, where $t\in[t_{k-1},t_k]$ and $1\leq k\leq n$. If $\alpha$ is not function of time then,
\begin{align*}
    \mathbb{P}\{\tau^{\alpha}_k \in ds\}=&\frac{\mathrm{d}}{\mathrm{d}t}F_{\tau_{k}^{\alpha}}(t)\\=& \frac{\mathrm{d}}{\mathrm{d}t}\Bigl \{\sum_{m=k}^{\infty}\frac{(-1)^m}{m!}\sum_{r=0}^{\infty}\frac{(-1)^r}{r!}\sum_{\substack{x_1,x_2,\ldots x_n \geq 0\\ x_1+x_2+\ldots+x_n=r}}\frac{r!}{x_1!x_2!\ldots x_n!}\lambda^{\sum_{i=1}^n\alpha_ix_i}(\alpha_1x_1+\ldots+\alpha_nx_n)_m\\& \prod_{k=1}^n (t_k-t_{k-1})^{x_k}\Bigl\}.
\end{align*}
\noindent
 The differentiation of $F_{\tau_{k}^{\alpha}}(t)$ with respect to $t_1, t_2,\ldots,t_k$, respectively and follows the uniform distribution can be write as
\begin{align}
   \frac{\mathrm{d}F_{\tau_{n}^{\alpha_n}}(t_n)}{\mathrm{d}t_n}=\begin{cases}
       h^{x_{k}+x_{k+1}-1}(x_{k+1}-x_{k})L & if~k=1,2,3\ldots,n-1.\\
       h^{x_{k}-1}x_{k}L & if~k=n.
   \end{cases}
\end{align}
where,  $L=\sum_{m=1}^{\infty}\frac{(-1)^m}{m!}\sum_{r=0}^{\infty}\frac{(-1)^r}{r!}\sum_{\substack{x_1,x_2,\ldots x_n \geq 0 \\x_1+x_2+\ldots+x_n=r}}\frac{r!}{x_1!x_2!\ldots x_n!}\lambda^{\sum_{i=1}^n\alpha_ix_i}(\alpha_1x_1+\ldots+\alpha_nx_n)_m  \prod_{k=1}^{n-2} (t_k-t_{k-1})^{x_k}$.

	\ifx
	Next we  evaluate \textit{pgf} of generalized space fractional Poisson process.
	\begin{theorem}
		Let $N(S^{\alpha(t)}(t)$ be the space fractional Poisson process
	\end{theorem}
	\begin{align*}
	\mathbb{E}[u^{N(S^{\alpha(t)}(t))}]=& \mathbb{E}[\mathbb{E}[u^{N(S^{\alpha(t)}(t))}|S^{\alpha(t)}(t)]] \\ =& \mathbb{E}[e^{-\lambda (1-u)S^{\alpha(t)}(t)}]\\ =& e^{-\int_0^t (\lambda (1-u))^{\alpha (s)}ds} \\ =& e^{-\sum_{k=1}^n \lambda^{\alpha_k}(1-u)^{\alpha_k}(t_k-t_{k-1})}\\ =& \prod_{k=1}^n e^{-\lambda^{\alpha_k}(1-u)^{\alpha_k}(t_k-t_{k-1})}
	\end{align*}
	which is the pgf of $ N(S^{\alpha_1}(t_1-t_0))+N(S^{\alpha_2}(t_2-t_2))+ \ldots N(S^{\alpha_n}(t_n-t_{n-1}))$\\
	Hence, \\
	$$ N(S^{\alpha(t)}(t) \stackrel{d}{=} N(S^{\alpha_1}(t_1-t_0))+N(S^{\alpha_2}(t_2-t_2))+ \ldots N(S^{\alpha_n}(t_n-t_{n-1})) $$
	\fi
	
\begin{theorem}
    The discrete Levy measure $\{\mathcal{V}_{\alpha(t)}(.)\}_{t\geq0}$ for \process\ $\gsfppvo$ is given by

    \begin{align}
        \mathcal{V}_{\alpha(t)}(.)=\sum_{m=0}^{\infty} \frac{(-1)^m}{m!}(\alpha_1x_1+\ldots+\alpha_nx_n)_m \delta_{\{m\}}(.)\int_0^{\infty} \exp{\{-\sum_{i=1}^{n}\lambda_{i}^{\alpha_{i}}(s_{i}-s_{i-1})\}} \nu(s)ds,
    \end{align}
    where $\delta_{\{m\}}(.)$ is Dirac measurement concentrated at m.
\end{theorem}
\begin{proof}
  The pmf of \process\ $\gsfppvo$ is
  \begin{align*}
   \mathbb{P}[\gsfppvo=m]=& \sum_{m=0}^{\infty} \frac{(-1)^m}{m!}\sum_{r=0}^{\infty}\frac{(-1)^r}{r!}\sum_{\substack{x_1,x_2,\ldots x_n \geq 0\\ x_1+x_2+\ldots+x_n=r}}\frac{r!}{x_1!x_2!\ldots x_n!}\lambda^{\sum_{i=1}^n\alpha_ix_i}(\alpha_1x_1+\ldots+\alpha_nx_n)_m \\& \prod_{k=1}^n (t_k-t_{k-1})^{x_k}.   
  \end{align*}
  To calculate Levy measures, we get
  \begin{align*}
      \mathcal{V}_{\alpha(t)}(.)=&\int_0^{\infty}\sum_{m=0}^{\infty}\mathbb{P}[\gsfppvo=m]\delta_{\{m\}}(.)\nu(s)ds\\=& 
      \int_0^{\infty}\sum_{m=0}^{\infty} \frac{(-1)^m}{m!}\sum_{r=0}^{\infty}\frac{(-1)^r}{r!}\sum_{\substack{x_1,x_2,\ldots x_n \geq 0\\ x_1+x_2+\ldots+x_n=r}}\frac{r!}{x_1!x_2!\ldots x_n!}\lambda^{\sum_{i=1}^n\alpha_ix_i}(\alpha_1x_1+\ldots+\alpha_nx_n)_m \\& \prod_{k=1}^n (s_k-s_{k-1})^{x_k} \delta_{\{m\}}(.)\nu(s)ds\\=&
      \sum_{m=0}^{\infty} \frac{(-1)^m}{m!}(\alpha_1x_1+\ldots+\alpha_nx_n)_m \delta_{\{m\}}(.)\int_0^{\infty}\sum_{r=0}^{\infty}\frac{(-1)^r}{r!}\sum_{\substack{x_1,x_2,\ldots x_n \geq 0\\ x_1+x_2+\ldots+x_n=r}}\frac{r!}{x_1!x_2!\ldots x_n!}\\&({\lambda_{1}^{\alpha_{1}}}(s_{1}-s_{0}))^{x_1}({\lambda_{2}^{\alpha_{2}}}(s_{2}-s_{1}))^{x_2}\ldots ({\lambda_{n}^{\alpha_{n}}}(s_{n}-s_{n-1}))^{x_n}\nu(s)ds\\=&
      \sum_{m=0}^{\infty} \frac{(-1)^m}{m!}(\alpha_1x_1+\ldots+\alpha_nx_n)_m \delta_{\{m\}}(.)\int_0^{\infty}\sum_{r=0}^{\infty}\frac{(-1)^r}{r!}\Bigl(\sum_{i=1}^{n}\lambda_{i}^{\alpha_{i}}(s_{i}-s_{i-1})\Bigl)^r \nu(s)ds\\=&
      \sum_{m=0}^{\infty} \frac{(-1)^m}{m!}(\alpha_1x_1+\ldots+\alpha_nx_n)_m \delta_{\{m\}}(.)\int_0^{\infty} \exp{\{-\sum_{i=1}^{n}\lambda_{i}^{\alpha_{i}}(s_{i}-s_{i-1})\}} \nu(s)ds.
  \end{align*}\qedhere
\end{proof}
\section{Aknowlagement}
We would like to acknowledge and thank Prof. Enzo Orsingher for proposing the central idea of this paper and for his careful review of the research paper. This article considered the variable-order of the stable subordinator as a function of time in a study of a generalized space-fractional Poisson process. His insightful suggestions substantially strengthened the quality of this work.

\bibliographystyle{abbrv}
\bibliography{researchbib}  
\end{document}